\documentclass[twocolumn,leqno]{article}
\usepackage[alphabetic,msc-links]{amsrefs}
  \pdfoutput=1

\usepackage[affil-it]{authblk}
\usepackage{graphicx}
\usepackage{color,bbm,subfig,caption,bm,fontenc,currvita}
\usepackage{amsmath,amsthm}
\usepackage{amsfonts}
\usepackage{amssymb,color,makeidx,enumerate,bbm}

\makeatletter
\renewcommand{\pod}[1]{\allowbreak\mathchoice
  {\if@display \mkern 18mu\else \mkern 8mu\fi (#1)}
  {\if@display \mkern 18mu\else \mkern 8mu\fi (#1)}
  {\mkern4mu(#1)}
  {\mkern4mu(#1)}
}

\def\x{{\bf x}}

\def\F{{\bf F}}

\def\bv{{\bf v}}

\def\f{{\bf f}}
\def\g{{\bf g}}
\def\h{{\bf h}}

\def\C{\mathbb{C}}
\def\1{\mathbbm{1}}

\def\bea{\begin{eqnarray}}
\def\eea{\end{eqnarray}}
\def\be{\begin{equation}}
\def\ee{\end{equation}}

\def\rar{\rightarrow}
\usepackage{layout}

\title{Operator is the Model}
\author{
 Igor Mezi\'c\footnote{
    Department of Mechanical Engineering and Department of Mathematics, University of California, Santa Barbara.
    }
  }

\begin{document}
\maketitle

Modeling of physical processes using dynamic evolution equations started in earnest with Isaac Newton. 
Ordinary differential equations (ODE's) came first and  Newton wrote to himself (in an anagram!),
that it was useful to solve them. Partial differential equations (PDE's) followed shortly thereafter, and these modeling methodologies still dominate 
applied mathematics. They utilize the concept of the dependent observable - 
that could be position and momentum in the case of Newton's gravitational models, or temperature in Fourier's
heat equation - and independent observables, such as time and spatial coordinate observables. 

Newton did his work on gravity with small amounts of data, relying on his (and perhaps Hooke's :-) brilliant intuition. In contrast, the end of the 20th and the beginning of the 21st century has seen a revolutionary 
increase in the availability of data. Indeed, we are in the middle of the {\it sensing} revolution,
where the word ``sensing" is used in the broadest meaning of data acquisition. Machine models (under the umbrella of Artificial Intelligence)
are used to analyze and make sense of this data, as we can witness from the current explosion of use of Large Language Models that rely on Deep Neural Networks technology, and more specifically on the idea of transformers. 
As these are typically vastly overparametrized (i.e. the number of weights in them is massive, in billions or even trillions - GPT-4 is rumored to have a trillion of those), an individual weight does not mean much. 
If humans are to correspond with machines intelligently, there is a need to extract models via which us humans can make our own sense of the data. 

However, often the assumptions that need to be satisfied for the underlying model to be an ODE or PDE, or another class of parsimonious (in the sense of deploying a small number of dependent and independent variables) are violated. 
To start with, there is observational noise, that could lead to another class of models - stochastic differential equations (SDE's)
provided some assumptions on the noise type are satisfied. But again, such assumptions might not be borne out. 

Koopman operator theory (KOT) has recently emerged as the main candidate for machine learning of {\it physics-based dynamical  processes }\cite{Mezic:2005,Budisicetal:2012,parmar2020survey,Mezic:2021}. I  propose here that its key paradigm is that  {\it ``the operator is the model"}. Namely, the assumption is that there exists a {\it linear} operator $U$, such that for {\it any}
observation $f$ of system dynamics $U$ enables prediction of the time evolution to the next observation $f^+$ using the equation
\be
f^+=Uf,
\ee
where $f$ is a function on some underlying state space $M$. As an example, the state space can be the space of position and linear momentum of a particle, and the observable function could be its energy. 
The operator $U$ is a property of the underlying dynamical process, and in that sense universal. However, it yields different outputs when applied to different observations - e.g. energy might be conserved over time while the position and the momentum are not. This change of setting - from dynamics on the state space, to dynamics on the space of observables $\cal{O}$ - led to a new modeling architecture that takes  $\cal{O}$ as its template.

Interestingly, as I show below, transformer architectures used in Large Language Models are in fact Koopman operator-based architectures. But in contrast to these, when applied with all the strength of the underlying theory, KOT  provides a powerful framework for unsupervised learning from small amounts of data, enabling self-supervised learning of {\it generative} models  that is much more in line with the theory of human learning than the machine learning methods of the second wave\footnote{DARPA classifies AI history into 3 waves \cite{prabhakar2017powerful} - roughly, the 1st wave is that of rule-based methodologies, the foundation of the 2nd  are supervised machine learning models and the 3rd wave is based on self-supervised, context-aware generative models.}. 

\section*{History}
\noindent Driven by success of the operator-based framework in quantum theory, Bernard Koopman proposed  \cite{Koopman:1931} to treat classical mechanics in a similar way, using the spectral properties of the composition operator associated with dynamical system evolution. Koopman  extended this study in a joint work with von Neumann in \cite{koopman1932dynamical}. Those works, restricted to Hamiltonian dynamical systems, did not attract much attention originally, as evidenced by the fact that between 1931 and 1990, the Koopman paper \cite{Koopman:1931} was cited 100 times, according to Google Scholar. This can be attributed largely to the major success of the geometric picture of dynamical systems theory in its state-space realization advocated by Poincar\'e. 
Out of today's 2000+ citations of Koopman's original work, \cite{Koopman:1931}, about 90\% come from the last 20 years. It was only in the 1990's  and 2000's that potential for wider applications of the Koopman operator-theoretic approach has been realized \cite{LasotaandMackey:1994,Mezic:1994,MezicandBanaszuk:2000,MezicandBanaszuk:2004,Mezic:2005,Rowleyetal:2009}. In the past decade the trend of applications of this approach has continued, as summarized in \cite{Budisicetal:2012,Mezic:2021}. This is partially due to the fact that strong connections have been made between the spectral properties of the Koopman operator for dissipative systems and the geometry of the state space. The Koopman operator framework is now in widespread use in machine learning - the number of papers in the field doubles every 5 years -
and its physical roots imbue it with {\it  interpretability}.

Even in the early work in \cite{Mezic:1994,MezicandBanaszuk:2000} and its continuation in \cite{MezicandBanaszuk:2004,Mezic:2005} there was an emphasis on utilizing the theory to find finite-dimensional models from data, as these concentrated on invariant subspaces of the operator - subspaces that are spanned by eigenfunctions. Finding an eigenfunction $\phi$ of the operator associated with a discrete-time, possibly  nonlinear process, yields a reduced model of the process, whose dynamics is governed by
\be
\phi^+=\lambda \phi,
\ee
and thus represents a  model of the dynamics\footnote{It is interesting to contrast transformation to eigenfunction space to the block of transformer architecture that nonlinearly transforms ``features" (in KOT language observables) into a latent space. For KOT, the latent space would be the space of eigenfunctions. The similarities between KOT and Large Language Model (LLM) architectures are further discussed below. }. The model yields {\it linear}
evolution of the function $\phi$. However, note that, in contrast with the standard setting in linear systems theory,  $\phi$ is a function on the underlying space $X$ that need not vanish on a linear manifold. Namely, the zero level set of $\phi$ can be a nonlinear set - e.g. a topological circle in the case of limit cycling dynamics \cite{mauroy2013isostables}. This is one example - where  a nonlinear attractor is described by a zero-level set of an eigenfunction -  of  the fact that  level sets of eigenfunctions on the original state space yield geometrically important objects such as invariant sets, isochrons and isostables \cite{mauroy2013isostables}. 
This led to realization that  geometrical properties can be {\it learned} from data, via computation of spectral objects, thus initiating a strong connection that is forming today between machine learning and dynamical systems communities \cite{li2017extended,yeung2019learning,lusch2018deep,takeishi2017learning}. The key notion driving these developments is that of representation of a -possibly nonlinear - dynamical system as a linear operator on a typically infinite-dimensional space of functions. This then leads to search of linear, finite-dimensional invariant subspaces, spanned by eigenfunctions. The idea is powerful: even multistable nonlinear systems can be represented this way \cite{Kvalheim:2023}. 

Numerous numerical methods were designed (e.g. \cite{MezicandBanaszuk:2004,Rowleyetal:2009,williamsetal:2015,das:2019} to find eigenfunctions and thus finite-dimensional models of the dynamics. But, it is of interest to invert the question and start {\it not} from the state-space model, but from the operator: $U$ is the property of the system - does it have a finite dimensional (linear or nonlinear) representation? In \cite{Mezic:2021}  the concept of {\it dynamical system representation} was formalized, enabling study of finite dimensional linear and nonlinear representations, learning, and the geometry of state space partitions. Instead of starting with the model, and constructing the operator, the finite-dimensional model is constructed {\it from} the operator. This enables  construction of models with {\it a-priori unknown physics}, used prominently in soft robotics \cite{bruderetal:2020,haggerty:2023}.
\section*{Operator Representations}
The modeling exercise typically starts with the catalogue of available observations. In Newton's case, these are positions and momenta of all the particles comprising the system. Let's label these observations on the abstract state of the system by  the vector $\f=(f_1,...,f_n)$, so we have $n$ different streams of data that we can organize into the $n\times m$ matrix $[\f(1),...,\f(m)]$, where $m$ is the number of  ``snapshots" of observations (in machine learning parlance observables are ``features"). For simplicity, we assume these snapshots are taken at regular time intervals and organized into columns sequentially. Note that 
\be
\f(k+1)=U\f(k),
\ee
where $U$ is the composition (Koopman) operator. We have made the assumption that the dynamics is evolving on some underlying state space $M$ (that we might not know) according to an unknown mapping $\x(k+1)=T(\x)$. The Koopman operator $U$ is then defined by
\be
 U\f=\f\circ T.
 \ee
  An interesting question to ask is: is there an $n\times n$ matrix $A$ such that
\be
U\f=A\f.
\ee
This is the case when $\f$ is in the span of $n$ (generalized) eigenfunctions of $U$ \cite{Mezic:2021}. More generally, we could ask: is there a map  $\F:\C^n\rar \C^n$ and observables (functions) $\g:X\rar \C^d$ such that
\be
U\g=\F(\g),
\label{repeig}
\ee
where $X$ is some latent space, and $\g=\g(\f)$. Typically, $d\geq n$.

\begin{figure}[ht]
\centering
\fbox{\includegraphics[  trim=0 150 70 100,clip=true,height=2.048in, width=3.1in
]{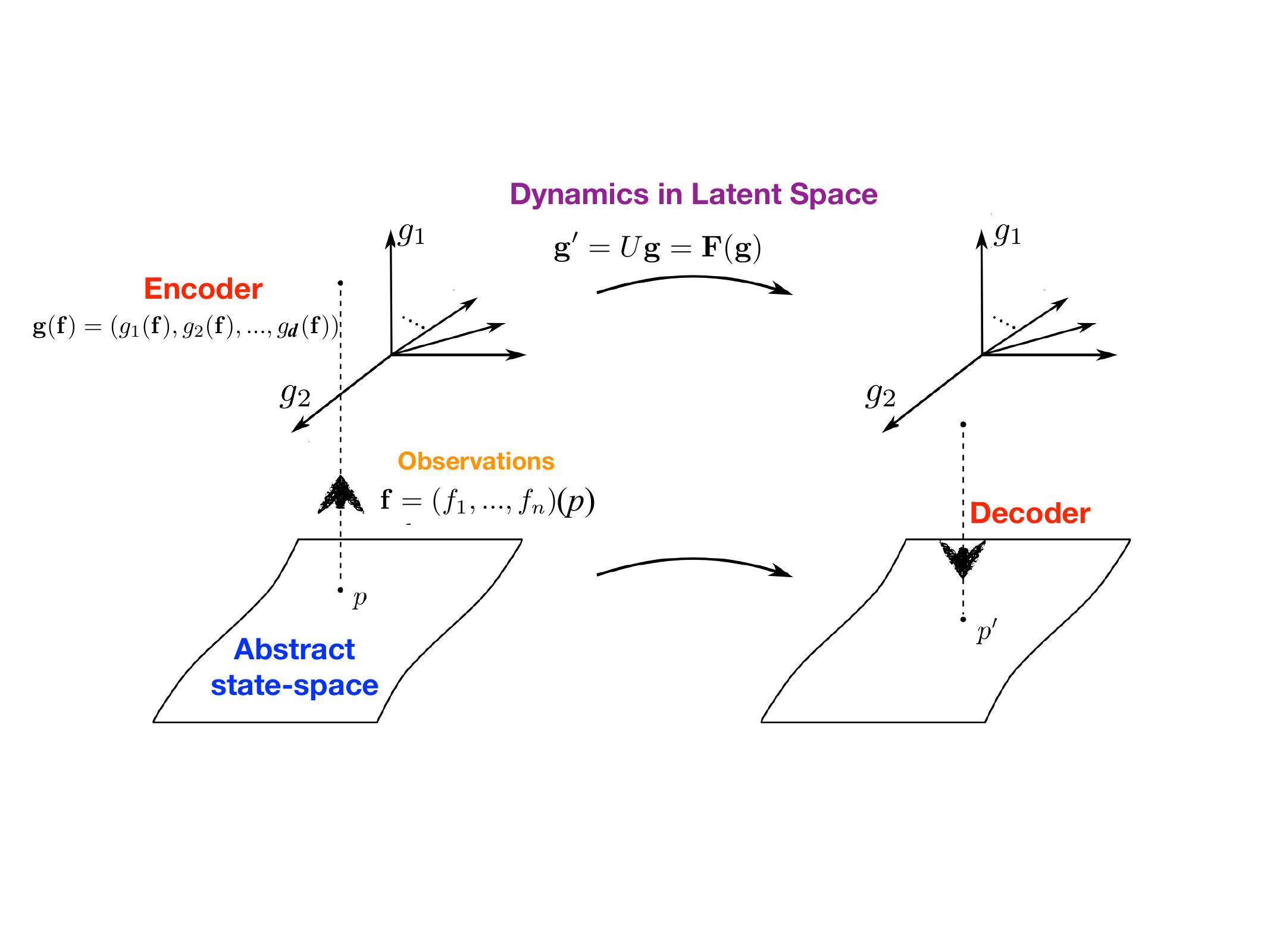}} \caption{The Koopman operator-based modeling architecture.}
\label{KoopEncDec}
\end{figure}

The modeling process is graphically represented in figure \ref{KoopEncDec}.

Taking our original observables $\f$ and setting $\g=\f$, it might be impossible to find such an $\F$. In that case, the observations do not provide us with a ``closure" i.e. we can not uniquely predict the next state of the observations from knowing the current state.

The problem of finding $(\F,\g)$ in (\ref{repeig}) was named the {\it representation eigenproblem} in \cite{Mezic:2021}. It reduces to the eigenvalue problem if we are seeking an eigenfunction
$\phi$ such that
\be
\phi(k+1)=\lambda\phi(k)
 \ee
 for some $\lambda \in \C$. In fact, there is  a precise result that tells us how the nature of the representation depends on the spectrum of the Koopman operator: finite linear representations are possible if the operator has discrete spectrum, while finite nonlinear representations are possible when the spectrum of the operator is continuous \cite{Mezic:2021}. 
 
 As an example of finding a nonlinear representation of the Koopman operator, consider discrete-time  dynamics given by $f_1=x,f_2=y$,
 \bea
 x^+&=&x+\sin x \\
 y^+&=&y+x,
 \eea
 (here we are pretending that the data is obtained by observing $f_1$ and $f_2$ but we do not know the underlying dynamics). We use $g_1(x)=x,g_2(x)=\sin x,g_3=y$. Now we have 
  \bea
 g_1^+&=&g_1+g_2 \label{eq:cl} \\
 g_2^+&=& g_2(g_1+g_2) \\
 g_3^+&=&g_3+g_1 
 \eea
The equation for $g_2$ does not have a linear combination of functions in the library on the right hand side.  To try to render this system of equations linear, without computing eigenfunctions,  we would need to start from introducing into the library the function
 \be
 g_4= g_2(g_1+g_2)=\sin(x+\sin x),
 \ee
  which would in turn 
 require more observables in order to close the system (i.e. in order for the left hand side depend linearly on the right hand side). On the other hand,
by just using $g_1$, we get a closed - albeit nonlinear - system, because the right hand side of (\ref{eq:cl}) can be written as $\F(\g)=F(g_1)=g_1+g_2(g_1)$. It turns out that the - infinite dimensional  - invariant subspace that consists of all the observables dependent on $f_1=x$ is the one on which the Koopman operator has a {\it nonlinear} representation $g_1' = g_1+g_2(g_1)$. In order to complete the model 
 we can set 
 \bea
 \g&=&(g_1,g_3), \\
  \F(\g)&=& (g_1+g_2(g_1),g_3+g_1).
  \eea
   If we were trying to approximate the action of the Koopman operator  on the data stream, and used $g_1=f_1,g_2=f_2=\sin (f_1),g_3=f_3$ as lifting to the latent space - thus utilizing two "state" functions $(f_1,f_3)$ and another function $f_2$ of the state function $f_1$ -
we could minimize
 \bea
 A&=&\min_{B}||G(k+1)-BG(k)|| \nonumber \\
 &=&\min_{B}||G^+-BG|| \nonumber 
 \eea
 where 
 \bea
 G=G(k)&=&[\g(1),\g(2),...,\g(k)],  \nonumber \\
 G^+=G(k+1)&=&[\g(2),\g(3),...,\g(k+1)] \nonumber
 \eea
  are  ``data" matrices whose columns are observables  evaluated at times $1,2,...,k$ and $2,3,...,k+1$, respectively. The solution, provided by
  \be
  A=G^+G^\dagger,
  \ee
  where $G^\dagger$ is the Moore-Penrose inverse of G,
yields (approximately)
 \be
 A\approx
 \left[\begin{matrix} 
 1 & 1& 0 \\
 * & * & * \\
 1 & 0 & 1 \\
 \end{matrix}\right].
 \ee
 Armed with the knowledge that $f_2=\sin f_1$, i.e. $f_2$ is functionally dependent on $f_1$, and that $ Uf_1$ is not a function of $f_3$, as indicated by the $0$ element of the first row of $A$, we conclude that there is a (reduced) 1-dimensional nonlinear representation of $U$ on the space of observables generated by $f_1$. From the last and first row of $A$  we conclude there is a (faithfull) 2-dimensional nonlinear representation. Thus, nonlinear representations can be extracted from data-driven computations.
 
 Note that here we have started from predefined observables $f_1=x,f_3=y$ and defined another observable $f_2=\sin f_1$ that we ``guessed" is important. Instead, we could have asked to minimize 
 \be
 (\beta^*,\gamma^*)=\min_{\beta,\gamma}||\g_\beta(k+1)-\F_\gamma(\g_\beta(k))||, 
 \ee
 where some - or all - of the components $g_j, F_k$'s are parametrized by neural networks with weights $\gamma,\beta$.  The interesting aspect of this is that it enables the concept of ``parenting" in learning. Namely, domain experts can suggest some of the key obsevables. A person that knows classical dynamics would suggest $\sin \theta$ as a good observable to be used for learning of the rigid pendulum dynamics, but use neural networks or time delay observables to learn appropriate observables for a soft pendulum \cite{bruderetal:2020,haggerty:2023} where physical laws are hard to derive. Thus a mixture of human-prescribed and machine-learned observables  is enabled.

 \section*{Koopman Operator Framework and Large Language Models}
 
 The current surge of attention in AI has been stimulated by the performance of Large Language Models underlying chatbots such as ChatGTP and Google Bard. In LLM's, language is considered as a dynamical system in which the next state of the system is determined by the previous states\footnote{This approach neglects the fundamental - {\it goal driven} - aspect of human intelligence - when we speak we typically do not just seek for the next word, but for the appropriate set of words to deliver the intended meaning. }. The embedding starts with  techniques such as One-Hot Encoding where each word is represented as a unit basis vector in the vector space whose dimension is that of a vocabulary. In the Koopman operator framework, these would be indicator observables on the set of words. The observables - features - that the LLM utilizes are combinations of time-delayed indicator observables. The transformer block then operates on the data matrix that contains features ordered in time (i.e. text order) and transforms the individual feature time sequence, followed by a transformation that nonlinearly combines transformed time sequences of features \cite{turner2023introduction}.
 
 This is exactly what the Koopman framework requires: abstract set elements are embedded into a Euclidian space\footnote{In the Koopman framework, observables can be complex - thus complex embeddings are enabled. LLM's use strictly real embeddings.}. Then, functions defined on that embedding are sought that can enable efficient prediction of  the dynamical system evolution. For example, time delayed observables can be used, also the common choice in LLM's. Filtering can be performed, producing linear combinations of such observables, equivalent to the first step in the transformer model \cite{turner2023introduction}. In the next step, a nonlinear transformation of these observables is sought, leading - if eigenvalues are found, in the case of discrete spectrum - to a linear representation, and to a nonlinear representation if the spectrum is continuous. In contrast to LLM's, Koopman operator based architectures are often computationally lean, because some of the transformations are not learned but hard-wired. Recent work discovered efficiency of the Koopman approach in language models \cite{orvieto2023resurrecting}. Interestingly, LLM-transformer architectures  were found inferior to    KOT-oriented architectures in time-series prediction \cite{sanchis2023easy,liu2023koopa}.
 \section*{Numerical Approximations}
Remarkably, the Koopman framework allows for finding of ``good" latent representations without ever knowing the operator. E.g., finding the eigenvalues and eigenfunctions associated with attractors can be pursued using harmonic analysis - and thus reduced to FFT's \cite{mezic:2022}
 or kernel methods \cite{das:2019}.
 Generalized Laplace Analysis methods extend these techniques to dissipative systems \cite{mezic:2022}. 
 
 A popular methodology for approximating the operator on a pre-determined feature set is the Extended Dynamic Mode Decomposition \cite{williamsetal:2015} that was recently incorporated into LAPACK codebase \cite{drmac2022lapack}. Machine learning methods have been developed to learn the features and maps at the same time \cite{li2017extended,yeung2019learning,lusch2018deep,takeishi2017learning}. All of these ultimately enable
 a reduced representation via the use of the spectral expansion \cite{Mezic:2005}, where the evolution of outputs (features) of interest $\h\in \C^h$ are represented (assuming discrete spectrum)
 by 
 \be
 \h(k+1)=\sum_{j=1}^N \lambda_j^k \phi_j(\x_0) \bv_j^{\h}
 \ee
 where $\lambda_j,\phi_j(\x_0)$ are a reduced set of eigenvalues and eigenfunctions and $\x_0$ is the initial condition that the data was initiated from. These are independent of the chosen set of features. The vector $\bv^{\h}_h$ is the Koopman mode obtained by projecting $\h$ onto the $j$-th eigenfunction, and depends on the selected  features (observables). The dynamics in the latent space of $\h$ is thus linear. The nonlinear dynamics of the original observations is then obtained via a (generally nonlinear) projection $\f=P_N(\h)$. The set $(\lambda_j,\phi_j,\bv^{\h}_j,j=1,...,N)$ is called the spectral triple, and can provide an extremely efficient 
 way of storing dynamic information - the state file needs to contain a reduced spectral triple set (size $(1+M+h)\times N$), where $M$ is number of initial conditions, and the projection map $P_N$.
 
 There is a recent powerful connection made to infinite-dimensional  numerical linear algebra methods \cite{colbrook:2023}. The resulting interplay between dynamical systems community and operator-theoretic numerical linear algebra community is guaranteed to broaden the horizons of both.
 \section*{Extensions and Relationship to Other Machine Learning Methods}
 
 Koopman based machine learning of dynamical models is particularly suitable for extension to control systems \cite{mezic:2004,mauroy2020koopman}. Koopman operators are defined that act on the tensor product structure of the lifting of state space and control space \cite{korda2018linear,cibulka2022dictionary}. Recent work in \cite{haseli2023modeling} has emphasized search for invariant subspaces in this tensor product space leading to models amenable to a plethora of control designs. Another extension is to machine learning of general nonlinear maps between different spaces \cite{Mezic:2021}. Stochastic effects have also been treated, as early as \cite{MezicandBanaszuk:2004,Mezic:2005} and later in \cite{takeishi2017subspace,vcrnjaric2020koopman,wanner2022robust,colbrook2023beyond}
 
 A number of connections have been made between ``pure" Koopman operator based methodologies, and other machine learning methodologies. The 
 version  of the framework with a pre-defined set of observables - starting from Schmid's DMD methods \cite{schmid:2010} - is conceptually equal to kernel methods popular in machine learning. The class of ARIMA models can be viewed as a subset of Koopman-based methods where only lifting in the input space has been performed, using time-delayed observables \cite{arbabi2017ergodic,mezic:2022}. As already mentioned, deep learning can be used to learn effective observables and connections to transformer architectures, widely used in LLM's have been made \cite{sanchis2023easy,geneva2022transformers}. 
 
 A popular methodology of reinforcement learning (RL) has been coupled to KOT modeling \cite{weissenbacher2022koopman}. However,  It is of interest to point out the fundamental difference of the approach to optimal control using KOT vs that used in RL: the exploration strategy in RL can lead to dangerous scenarios. In KOT approach to optimal control, the model is formed first, assuring that only safe scenarios are executable. Then a cost function is specified, enabling optimization of the task while safety is preserved.
  
 Because of the explicit treatment of the time-dimension, Koopman operator models - that fit the nomenclature of foundational models in generative AI - are well suited for dealing with causal inference in the sense of Pearl \cite{pearl2019seven}. For example, counterfactual questions such as ``What if I
had acted differently?" can be answered using a Koopman control model. In fact, all of Pearl's obstacles to developing autonomous systems
that exhibit human-level intelligence - robustness, adaptability, explainability (interpretability) and cause-effect relationships can be resolved using generative Koopman control models. The applications of the methodology are extensive - starting from fluid mechanics \cite{Mezic:2005,Rowleyetal:2009,mezic2010new}, continuing with power grid \cite{susuki2016applied}, the methodology has now penetrated most fields where dynamics is important, including recent advances in synthetic biology \cite{yeung2019learning} and soft robotics \cite{bruderetal:2020,haggerty:2023}. It was even used to model the Starcraft game \cite{avila2021game}! This is partly due to the effectiveness of developed machine learning algorithms, but also due to the depth of the underlying theory that enhances interpretability, prevalent in applied mathematics, but missing in much of the  ``pure" machine learning approaches.

 Despite all of the described progress, there is still much to do, and the current decade is going to be an exciting one for this growing set of  data-driven AI methodologies for discovering models of dynamical processes.
\bibliography{KvN,MOTDyS,MOTDyS2}

\end{document}